\documentclass{cpc}

\usepackage{amssymb,amsmath}

\newcommand{\rest}{\upharpoonright}

\newcommand{\CST}{Central Sets Theorem}
\newcommand{\FS}{{FS}}       
\newcommand{\FP}{{FP}}       
\newcommand{\SCC}{Stone-{\v C}ech compactification}

\renewcommand{\P}{\mathcal{P}}         
\newcommand{\N}{\mathbb{N}}             
\newcommand{\Z}{\mathbb{Z}}             

\renewcommand{\Cup}{\bigcup}
\renewcommand{\Cap}{\bigcap}
\renewcommand{\subset}{\subseteq}

\newcommand{\app}{^{\smallfrown}}
\newcommand{\dom}{{dom\ }}

\newcommand{\kw}{{<\omega}}
\newcommand{\me}{{-1}}
\newcommand{\fin}{\P_f (\omega)}
\newcommand{\seq}[1]{{\langle {#1}_{n}\rangle_{n=0}^\infty}}

\renewcommand{\Cup}{\bigcup}
\renewcommand{\Cap}{\bigcap}
\newcommand{\nhat }[1]{\{1,2,\ldots,#1\}}
\newcommand{\ohat }[1]{\{0,1,\ldots,#1\}}

\begin{document}

\title{A Multidimensional Central Sets Theorem}
\author[M. Beiglb\"ock]{\spreadout{M. BEIGLBOECK}
\thanks{The author thanks the Austrian
Science Foundation FWF for its support through
Project no. S8312. and Project no. P17627-N12}
\\
       \affilskip Department of Discrete Mathematics and Geometry, 
        TU Vienna,\\       \affilskip 1040 Vienna, Austria,
       mathias.beiglboeck@tuwien.ac.at    }

\maketitle 
\abstract{ 
The Theorems of Hindman and van der Waerden belong to the classical 
theorems of partition Ramsey Theory. The Central Sets Theorem is a 
strong simultaneous extension of both theorems that  
applies to general commutative  
semigroups. We give a common extension of the Central Sets Theorem 
and 
Ramsey's Theorem.  
} 
 
\section{Introduction}  
 
Van der Waerden's Theorem (\cite{vdW}) 
 states that for any partition of the positive integers 
$\N$ one of the cells of the partition contains arbitrarily long arithmetic  
progressions. 
 
To formulate Hindman's Theorem (\cite{H}) and the \CST\ we set up some notation. 
By  $\fin$ we denote the set of all finite nonempty subsets of  
$\omega= \N \cup \{0\}.$ 
 For a sequence $\seq x$ in $\N$  
we put \emph{$ \FS (\seq x):= 
\{ \sum_{t\in\alpha} x_t: \alpha \in \fin\}.$} 
A set $A\subset \N $ is called an \emph{IP-set} iff there exists a sequence 
$\seq x $ in $\N$ such that $ \FS (\seq x)\subset A$. (This 
definitions make perfect sense in any semigroup $(S,\cdot)$ 
 and we indeed plan to  
use them  in this context. $\FS$ is an abbriviation of \emph{finite sums} 
and will be replaced by $\FP$ if we use multiplicative  
notation for the semigroup operation.) 
Now Hindman's Theorem  states that in any finite  
partition of $\N$ one of the cells 
is an IP-set.      
 
K. Milliken and A. Taylor (\cite{M,T}) 
found a quite natural common extension of 
the Theorems of Hindman and Ramsey:  
For a sequence $\seq x$ in $\N$ 
and $k\geq 1$ put  
$ [\FS (\seq x)]_<^k:= 
\left\{ \left\{\sum_{t\in\alpha_1} x_t,\ldots, 
 \sum_{t\in\alpha_k} x_t\right\}  :  
\alpha_1 <  \ldots <\alpha_k \in \fin\right\},$ 
where 
we write $\alpha < \beta$ for $\alpha, \beta\in \fin$ iff  
$\max \alpha < \min \beta$. For an arbitrary set $S$ let $[S]^k$ be the 
 set of 
all finite subsets of $S$ consisting of exactly $k$ elements. 
If $[\N]^k = \Cup_{i=1}^r A_i $ then there exist 
 $i\in \nhat r$ and  
a sequence $\seq x$ in $\N$ such that  
$\FS [\seq x]_<^k\subset A_i.$  
 
Let $\Phi$ be the set 
 of all functions $f:\omega \rightarrow \omega$ such that  
$f(n) \leq n$ for all $n\in \omega$. Then our main theorem may be 
 stated as follows:  
 
\begin{theorem}\label{meinmaintheorem} 
Let $(S,\cdot)$ be a commutative semigroup  
and assume that there exists a non principal 
minimal idempotent in $\beta S$. For each $l\in \N$, let  
$\langle y_{l,n}\rangle_{n= 0}^\infty $ be a sequence in $S$. 
Let $k,r \geq 1$ and let $ [S]^k = \Cup_{i=1}^r A_i$. There exist  
$i\in \nhat r$, a sequence $ \seq a $ in $S$ and a sequence  
$\alpha_0<\alpha_1< \ldots$ in $ \fin$   such that for each $g\in \Phi$,  
$ \left[ \FP \left( \left\langle a_n  \prod_{t\in \alpha_n}  
   y_{g(n),t}\right \rangle_{n= 0}^\infty\right)\right]^k_< \subset  A_i. $ 
\end{theorem} 
We will review some properties of the \SCC\ as well as 
the definition of  a minimal idempotent in the next chapter.  
In the case $k=1$ the somewhat odd assumption 
that $\beta S$ should contain a non 
 principal minimal idempotent is not needed. 
In general this condition will be satisfied  
 if $S$ is  \emph{weakly (left)  
cancellative}, i.e. for all $u,v \in S$ the set $\{s\in S : us=v\}$  
is finite 
and $S$ itself is infinite 
(see \cite{HS}, Theorem 4.3.7). In particular the conclusion of 
Theorem \ref{meinmaintheorem} holds  in the semigroups $(\N, + ), (\N, \cdot), 
(\fin, \cup)$.  
 
The case $k=1$ of Theorem \ref{meinmaintheorem} is exactly the \CST.  
(More precisely this is the version stated  
in \cite{HS}, Corollary 14.12. A  
discussion on the origin of the \CST\ can also be found there.) 
By further specifying $(S,\cdot)= (\N,+)$ and 
 $ \langle y_{l,n}\rangle_{n= 0}^\infty=\langle l, l,\ldots\rangle$ 
we get that all finite sums of elements of the 
arithmetic progressions $a_n, a_n + |\alpha_n|, \ldots, 
 a_n + n |\alpha_n|, n\geq 0$   
are guaranteed to be monochrome. 
   
 Theorem   \ref{meinmaintheorem} may be 
 seen as a generalization of the \CST\ 
in the same sense  as the Milliken-Taylor Theorem  
is a multidimensional version of Hindman's Theorem.

\section{Preliminaries on ultrafilters} 
 
For a set $S$ let $\beta S$ be the set of all ultrafilters on $S$. For  
$s\in S$ we will identify $s$ with the principal ultrafilter of all  
subsets of $S$ that  
contain $s$.  If $(S,\cdot)$ is a semigroup, the operation $.\cdot.$ 
may be extended to a semigroup operation on $ \beta S$ by defining 
\begin{equation}\label{multdef} 
  A \in p \cdot q  \ :\Leftrightarrow \ \{s\in S : s^\me A \in q\}\in p.       
\footnote{$S$ is a semigroup, so $s$ might not have an inverse. We may 
avoid this obstacle by defining $s^\me A:= \{t\in S: st\in A\}.$} 
\end{equation} 
If $\beta S$ is properly topologized it turns out to be the \SCC\ of $S$  
(where we regard $S$ to be a discrete space).  
It can be shown that the operation  
$ .\cdot.: \beta S \times \beta S \rightarrow  
\beta S$  defined in (\ref{multdef}) is the unique   
extension of $.\cdot. :S\times S  
\rightarrow S$, such that for each $s\in S$ and  
each $q\in \beta S$  the functions  
$\lambda_s,\rho_q:\beta S \rightarrow \beta S$ defined by 
 $ \lambda_s (r) := sr, \rho_q (r) := rq $ are continuous.   
 
Applications of the algebraic structure of $\beta S$ in partition Ramsey Theory  
are abundant. Examples are simple proofs of the theorems of Hindman and  
van der Waerden: 
 
Idempotent ultrafilters (i.e. ultrafilters $e\in \beta S$ 
satisfying $e e =e$) turn out to be tightly connected with  IP-sets  
in $S$: A subset $A$ of $S$ is an IP-set iff there is an idempotent $e\in \beta  
S$ such that $A\in e$.  
By a theorem of Ellis $\beta S$ always 
 contains an idempotent  ultrafilter 
$e$ and by the ultrafilter properties of $e$ for any partition  
$A_1, A_2, \ldots, A_{r}$ of $S$ there exists an $i\in\nhat r$ such  
that $A_i\in e$. Thus  
$A_i$ is an IP-set.  
 
$\beta S$ always has a smallest (two-sided) ideal which will be denoted by  
$K(\beta S)$. It turns out that for $(S,\cdot)=(\N, +)$ the elements of 
$K(\beta \N)$  are well suited for van der Waerden's Theorem.  
 
Idempotents  in $K(\beta S)$ (which are always present) are called 
\emph{minimal idempotents}. 
Not at all surprisingly minimal idempotents are   
particularly interesting for combinatorial applications. Subsets of 
$S$ which are contained in some minimal idempotent are called \emph{central sets} 
and that these sets satisfy the conclusion of the \CST\ reveals the source of 
the theorem's name. 
 
See \cite{HS} for an elementary introduction to the semigroup $\beta S$ as 
well as for the combinatorial applications mentioned in this section. 
 
If $S$ is an infinite set an arbitrary non principal ultrafilter 
$p\in \beta S$ may be used to give a proof of Ramsey's Theorem.  
(This proof  is by now classical. See \cite{C} p.39 
for a discussion of its origins.) It's an idea of V. Bergelson and  
N. Hindman that in the case $S=\N$, something 
might be gained by using an ultrafilter with special  algebraic  
properties. Via this  approach in \cite{BH} a short proof of 
the Milliken-Taylor Theorem is given and  
a very strong simultaneous generalization 
of Ramsey's Theorem and  numerous single-dimensional Ramsey-type Theorems 
(including van der Waerden's Theorem) is obtained.   
Our proof is a variation on this idea.

\section{The proof of the main theorem} 
 
The following Lemma is the basic tool in the ultrafilter proof of  
Ramsey's theorem: 
\begin{lemma}\label{basiclemma} 
Let $S$ be a set, let $e\in \beta S \setminus S$, 
 let $ k,r \geq 1 $, and let $[S]^k = \Cup_{i=1}^r A_i$. For each  
$i\in \nhat r$, each $t\in \nhat k$   
and each $E\in [S]^{t-1}$, define $B_t(E,i)$ by downward induction  
on $t$: 
\begin{enumerate} 
\item[(1)] For $E\in [S]^{k-1}, B_k (E,i)  
  := \{ y\in S\setminus E: E\cup \{ y\} \in A_i \}. $ 
\item[(2)] For $1\leq t <k$ and $E\in [S]^{t-1},$ 
$ B_t(E,i) := \{ y\in S\setminus E : B_{t+1} (E\cup \{y\},i)\in e\}.$ 
\end{enumerate}  
Then there exists some $i\in \nhat r$  
such that $B_1(\emptyset, i)\in e$. 
\end{lemma} 
\begin{proof} 
 For each $E\in [S]^{k-1}$ one has  
$S= E\cup \Cup_{i=1}^r B_k(E,i)$, so there  
exists $i\in \nhat r$ such that $ B_k (E,i) \in e$.  
Next let $E\in [S]^{k-2}$ and $y\in S\setminus E$. Then there 
exists $i\in \nhat r$ such that $B_k (E\cup \{y\},i) \in e$.  
Thus $S= E\cup \Cup_{i=1}^r B_{k-1} (E,i)$.  
\\ After iterating this argument $k-1$ times we achieve  
  $S = \emptyset \cup \Cup_{i=1}^r B_1(\emptyset,i)$ which clearly proves  
the statement.  
\end{proof} 
To formulate our key lemma we need to introduce some notation: 
Let  $S$ be a set and  put $S^\kw = \Cup_{n= 0}^\infty S^{\{0,\ldots, n-1\}}$. 
A non empty set $T\subset S^{\kw}$ is a \emph{tree} in $S$ 
  iff for all $f\in S^\kw$, $g\in T$ such that 
$\dom f\subset \dom g$, $g_{\rest \dom f} =f$ one has $f\in T$.   
 We will identify a function $ 
f\in S^{\ohat {n-1}}$ with the sequence $\langle  
f(0),f(1), \ldots, f(n-1)\rangle.$  
 If $s\in S$ then $f\app s :=\langle f(0),f(1),\ldots,f(n-1), s\rangle.$ 
 For $f\in S^\kw$ we put $T(f) := \{ s\in S : f\app s \in 
  T\}$.  
 
\begin{lemma}\label{baumramsey} 
Let $(S,\cdot)$ be a semigroup such  
that there exists an idempotent $e\in \beta S \setminus S$, 
 let $k,r\geq 1$ and 
 let $[S]^k =\Cup_{i=1}^r A_i $.  
 Then there exist  $i \in \nhat r$ and a {tree} $T\subset S^\kw$ such  
that for all $f\in T$ and $\alpha_1 <\alpha_2< \ldots < \alpha_{k} \subset \dom 
f, \alpha_i\in \fin$ one has: 
\begin{enumerate} 
 \item[(1)] $T(f)\in e$. 
\item[(2)]  $ \left\{ \prod_{t\in \alpha_1} f(t), 
\prod_{t\in \alpha_2} f(t), \ldots , 
 \prod_{t\in \alpha_k} f(t)\right\}  \in A_i.$ 
\end{enumerate} 
\end{lemma}

In the proof we will employ some basic properties of 
idempotent ultrafilters. 
For a set $A\subset S$ we have 
$A\in e= e e$ iff $\{s\in S:s^\me A\in e\}\in e$ 
by the definition of the multiplication in $\beta S$. 
For $A\in e$ let  
$$A^\star:=\{s\in A:s^\me A\in e\}= A\cap \{s\in S:s^\me A\in e\}\in e$$ 
Then $t^\me A^\star\in e $ for all $t\in A^\star$: 
\\ It is clear that for $t\in A^\star, t^\me A\in e$. Furthermore  
$$t^\me\{s\in S: s^\me A\in e\}=\{s\in S: s^\me(t^\me A )\in e\} \in e.$$ 
Thus 
in fact $t^\me A^\star=t^\me A\cap t^\me\{s\in S: s^\me A\in e\}\in e $. (This is \cite{HS}, Lemma 4.14.)

\begin{proof} 
Let $i\in \nhat r$ be such that $B_1(\emptyset,i)\in e$. (We use the  
notation of Lemma \ref{basiclemma}. 
Since $i$ will be fixed in the rest of the proof, we will  
suppress it and write $B_r(E)$ instead of $B_r(E,i)$.) 
We will inductively  
construct an increasing sequence of trees  
$\langle T_n\rangle_{n=0}^\infty$, satisfying for 
each $n\geq 0$,  
$T_n=\{f_{\rest \nhat{n-1}}: f\in T_{n+1}\}$ 
such that the for each $f\in T_n$ the following holds: 
\begin{enumerate} 
\item[(i)] 
If $\dom f \subset \ohat {n-2}$ then $T_n(f)\in e$. 
\item[(ii)] 
If $\alpha_1,\alpha_2, \ldots,\alpha_r\in \fin,r\in \nhat k$ 
satisfy  
$\alpha_1< \alpha_2<\ldots< \alpha_r\subset \dom f$ 
and if 
$x_i=\prod_{t\in \alpha_i} f(t)$ then 
$x_r \in B_r(\{x_1, x_2,\ldots,x_{r-1}\})^\star$\footnote{For 
$r=1$ this is meant to be $B_1(\emptyset)^\star$.}.  
\end{enumerate} 
Trivially we may put $T_0= \{\emptyset\}$. Assume now that  
$T_0, T_1,\ldots, T_n$ have already been defined.  
Fix $f\in T_n$ with $\dom f =\ohat{n-1}$. For 
$\alpha_1<\alpha_2<\ldots<\alpha_r\subset \dom f$ let  
$x_i=\prod_{t\in \alpha_i} f(t)$. By assumption  
$x_r \in B_r(\{x_1,x_2,\ldots,x_{r-1})\}$ and thus  
$B_{r+1} (\{x_1,x_2,\ldots, x_r\})\in e$ for $r\in \nhat {k-1}$. 
Since   $x_r\in B_r(\{x_1,x_2,\ldots, x_{r-1}\})^\star$ we have 
 $x_r^\me B_{r}(\{x_1,x_2,\ldots,x_{r-1}\})^\star\in e$  
for $r\in \ohat k$. Define $T_n(f)$ to be the  
intersection of all sets $B_{r+1} (\{x_1,x_2,\ldots, x_r\})^\star,$  
$r\in \ohat {k-1}$ 
and $x_r^\me B_{r}(\{x_1,x_2,\ldots,x_{r-1}\})^\star 
, r\in \ohat k$ such that indeed $T_n(f)\in e$. Using this 
put $T_{n+1}= T_n\cup \{f\app t:f\in T_n, \dom f= \ohat {n-1}, t\in T_n(f)\}$. 
It is not hard to verify that this implies that the inductive construction  
can be continued: This is only interesting for $\dom f=\ohat n$  
and $n\in \alpha_r$ (where  $r\in \nhat k$).  
Fix $f':\ohat {n-1}\rightarrow S$ such that 
${f'}\app f(n)=f$. 
If $\alpha_r=\{n\}$, 
$x_r= f(n)\in T_n(f') \subset B(\{x_1,x_2\ldots,x_{r-1}\})^\star$ 
so we are done. 
If $\alpha_r=\alpha_r'\cup \{n\}$ for some non empty 
$\alpha_r'\subset \ohat {n-1}$ we have 
$f(n)\in T_n(f')\subset \left(\prod_{t\in \alpha_r'}f'(t)\right)^\me  
B_r(\{x_1,x_2,\ldots,x_{r-1}\})^\star$ and this implies 
$x_r =\prod_{t\in \alpha_r} f(t)  
\in B_r(\{x_1,x_2,\ldots,x_{r-1}\})^\star.$ 
 
Finally put $T=\Cup_{n=0}^\infty T_n$. Obviously $T(f)\in e$ for  
all $f\in T$. Since $\prod_{t\in \alpha_k} f(t)\in B_k\left(\left\{  
\prod_{t\in \alpha_1} f(t), \ldots, \prod_{t\in \alpha_{k-1}}f(t) 
\right\} \right)$ for  
all $f\in T$ and $\alpha_1<\alpha_2<\ldots<\alpha_k\subset \dom f$ 
we see that (2) holds. 
\end{proof} 

From this Lemma one may directly derive the following strong 
 version of the Milliken-Taylor Theorem: 
\begin{corollary}\label{meinmt} Let $k,r\geq 1$, 
 let $(S,\cdot)$ be a semigroup,  
let $\seq x $ be a sequence in $S$ and  
let $[S]^k=\Cup_{i=1}^r A_i$.  
Assume  that for every idempotent $s\in S$ there exists some  
$m\in \N$ such that $s\notin  
\FP (\langle x_n\rangle_{n= m}^\infty)$. 
  Then there exist $i\in \nhat r$  
and a sequence $\alpha_0< \alpha_1< \ldots$ in $\fin$ such that 
$\left[\FP \left( \left\langle  
      \prod_{t\in \alpha_n} x_t \right\rangle_{n= 0}^\infty 
               \right)\right]^k_<\subset A_i.$   
\end{corollary} 
\begin{proof} 
By \cite{HS}, Lemma 5.11 there exists an idempotent $e\in \beta S$, 
 such that 
for all $m\geq 0$, $\FP (\langle x_n\rangle_{n= m}^\infty )\in e$ and  
by our assumption we have  
$e\in \beta S \setminus S$. Let $i\in \nhat r$ and 
 $T\subset S^\kw$ be as provided  
by Lemma \ref {baumramsey}.  
 We have $T(\emptyset) \cap \FP (\seq x)\in e$. In particular 
this set is not empty, so we may choose $\alpha_0\in \fin$ such that 
$\prod_{t\in \alpha_0} x_t \in T(\emptyset)$. Let $ m_0:= \max \alpha_0$.  
As before $T\left(\left( \prod_{t\in \alpha_0} x_t\right) \right) \cap  
 \FP \left( \langle x_n\rangle_{n=m_0}^\infty\right) \in e$, so we find  
$ \alpha_1 > \alpha_0, \alpha_1\in \fin$ such that  
$ \prod_{t\in \alpha_1}x_t \in T\left(  
   \left\langle\prod_{t\in \alpha_0} x_t \right\rangle \right)$.    
By continuing in this fashion we achieve a sequence  
with the required properties.  
\end{proof}  
We remark that our restriciton on the idempotents  
contained in $\FP(\seq x)$ cannot be dropped: 
Consider for example $(S,\cdot)= (\Z,+)$ 
 and $\seq x= \langle 0,0, \ldots\rangle$: 
In this case $[\FP (\seq x)]_<^k=\{\{0\}\}$ for 
any $k\in \N$. 
 
Another possibility to avoid this difficulty is presented  
in 
\cite{HS}, 
 Corollary 18.9: 
Instead of partitions of 
$[S]^k$, partitions of $\Cup_{i=1}^k[S]^i$ are considered there. 
 
In the proof of Theorem \ref{meinmaintheorem} we will require the following: 
\begin{theorem}\label{lemmm2} 
Let $(S,\cdot)$ be a commutative  semigroup, let 
$A\in e\in K(\beta S)$,  
let $l\in \N$ and for each $j\in \ohat {l-1}$ let  
$\langle y_{j,n}\rangle_{n= 0}^\infty $ be a sequence in $S$. 
 Then there exist  
 $a \in S$ and   $\alpha \in \fin$    such that  
$a  \prod_{t\in \alpha} y_{j,t} \in A$ for each $j\in \ohat{l-1}$.  
\end{theorem}  
Theorem \ref{lemmm2} is a special case of the \CST\  and may easily  
be derived from the Hales-Jewett Theorem (\cite{HJ}). 
 
We are now able to prove our main Theorem: 
\begin{proof}[Proof of Theorem \ref{meinmaintheorem}.] 
Fix a minimal idempotent $e\in \beta S \setminus S$.  
Let $i\in \nhat r$ and $T\subset S^\kw$ be  
as provided by lemma \ref{baumramsey}. 
We will inductively construct  sequences  
$\seq a$ in $S$ and  $\alpha_0 < \alpha_1 < \ldots$ in $\fin$  
 such that for all $n\in 
\N$ and all $g\in \Phi$: 
\begin{equation}\label{cstindvor} 
 \mbox{$ \left\langle a_0 \prod_{t\in \alpha_0 }y_{g(0),t},  
a_1 \prod_{t\in \alpha_1 }y_{g(1),t}, 
\ldots , a_{n-1}  \prod_{t\in \alpha_{n-1}} y_{g(n-1),t} 
  \right\rangle \in T. $} 
\end{equation} 
By the properties of $T$ this is sufficient to proof the Theorem.  
 
Assume that  $a_0 ,a_1, \ldots , a_{n-1} \in S$ und $\alpha_0< \ldots  < 
\alpha_{n-1} \in \fin $ have already been constructed such that (\ref{cstindvor}) 
is true for all $g\in \Phi$.  
We have 
$$G_n:= \Cap_{g\in \Phi}
T\left( \left\langle a_0\!\! \prod_{t\in \alpha_0}\!\!y_{g(0),t},\  
a_1\!\! \prod_{t\in \alpha_1}\!\!y_{g(1),t}, 
\ldots , a_{n-1}\!\!\!\!\!\! \prod_{t\in \alpha_{n-1}}\!\!\!\!y_{g(n-1),t} 
  \right\rangle\right)\in e. $$   
Let $m:=\max \alpha_{n-1}$. 
   By applying  Theorem \ref{lemmm2} to the set  
   $G_n$ and the sequences $\langle y_{0,k}\rangle_{k=m}^\infty, 
\langle y_{1,k}\rangle_{k=m}^\infty, \ldots 
  ,\langle y_{n,k}\rangle_{k=m}^\infty$ we find $a_n \in S$ und $ 
  \alpha_n \in \P_f(\N)$, $\alpha_n > \alpha_{n-1}$ such that 
$  a_n  \prod_{t\in \alpha_n} y_{0,t} , 
a_n  \prod_{t\in \alpha_n} y_{1,t} , 
 \ldots , a_n \prod_{t\in \alpha_n} 
  y_{n,t} \in  G_n. $  
\\ Thus for all $g\in \Phi$,  
 $ \left\langle a_0 \prod_{t\in \alpha_0}y_{g(0),k}, 
a_1 \prod_{t\in \alpha_1}y_{g(1),k}, \ldots  
, a_{n} \prod_{t\in \alpha_n}y_{g(n),k} 
  \right\rangle \in T, 
$ as we wanted to show. 
\end{proof} 
 
We conclude this section by giving a strengthening of Theorem  
\ref{meinmaintheorem} that applies to partitions of 
the spaces $[S]^1,[S]^2,\ldots, [S]^k$ simultaenously. It 
is not hard to verify that a similar extension of Corollary \ref{meinmt} 
is also valid. To avoid confusion about indices we use colourings 
instead of partitions.  
 
\begin{corollary}\label{meinmaintheorem+} 
Let $(S,\cdot)$ be a commutative semigroup  
and assume that there exists a non principal 
minimal idempotent in $\beta S$. For each $l\in \N$, let  
$\langle y_{l,n}\rangle_{n= 0}^\infty $ be a sequence in $S$. 
Let $k \geq 1$ and assume that  
for each $m\in \nhat k$,  $ [S]^m$ is finitely coloured. 
 There exist  
 a sequence $ \seq a$ in $S$, a sequence  
$\alpha_0<\alpha_1< \ldots$ in $ \fin$ and 
for each $m\in \nhat k$ a monochrome set $A^{(m)}$ 
   such that for each $g\in \Phi$ and each $m\in \nhat k$,  
$ \left[ \FP \left( \left\langle a_n \prod_{t\in \alpha_n}  
   y_{g(n),t}\right\rangle_{n= 0}^\infty\right)\right]^m_< \subset  A^{(m)}. $ 
\end{corollary} 
 
\begin{proof} 
We describe two ways to prove  
Corollary \ref{meinmaintheorem+}: 
\\  
Fix a linear ordering $\prec$ on $S$. 
For $m\in \nhat k$ let 
$ f^{(m)}:[S]^m\rightarrow \nhat {r_m} $ be the colouring at hand.  
Define a colouring 
$g^{(m)}:[S]^k \rightarrow \nhat {r_m}$ by 
 letting $g^{(m)}(\{x_1,x_2,\ldots, x_k\})= 
f{(m)} (\{x_1,x_2,\ldots,x_m\})$, where  
$\{x_1,x_2, \ldots, x_m \}$ are the $m$ smallest elements of  
$\{x_1,x_2,\ldots,x_k\}$ with respect to $\prec$.  
Then apply Theorem \ref{meinmaintheorem} to the 
colouring  
\begin{eqnarray*} 
f:[S]^k&\rightarrow& \nhat {r_1} \times \nhat {r_2} 
 \times \ldots \times \nhat r_k\\ E&\mapsto& (g^{(1)}(E),g^{(2)}(E),\ldots, 
g^{(k)}(E)).\end{eqnarray*} 
It is clear that the resulting sequences  
$\seq a$ and $\seq \alpha $ satisfy the conclusion of  
Corollary \ref{meinmaintheorem+}. 
 
The more complicated way to prove Corollary \ref{meinmaintheorem+} 
is to start by extending Lemma \ref{baumramsey}.  
Pick a minimal idempotent $e\in \beta S \setminus S$.  
Choose by Lemma \ref{baumramsey} 
 for each $m\in \nhat k$ a monochrome  
set $ A^{(m)} \subset [S]^m $ and a tree $T^{(m)}\subset S^\kw$ such that  
for all $f\in T^{(m)}$ and all $\alpha_1<\alpha_2<\ldots<\alpha_m\subset 
\dom f, \alpha_i\in \fin$ one has $T^{(m)}(f) \in e $ and 
$\{ \prod_{t\in \alpha_1} f(t), \prod_{t\in \alpha_2}f(t), 
\ldots, \prod_{t\in \alpha_m} f(t)\}\in A^{(m)}.$ But then 
$T:= \Cap_{m=1}^k T^{(m)}$ is a tree such that for all $f\in T, 
T(f)\in e$ and for all $ m\in \nhat k$ and all  
$\alpha_1<\alpha_2<\ldots<\alpha_m \subset \dom f, \alpha_i\in \fin $,  
 $\{ \prod_{t\in \alpha_1} f(t), \prod_{t\in \alpha_2}f(t), 
\ldots, \prod_{t\in \alpha_m} f(t)\}\in A^{(m)}.$ By  
 performing the proof of Theorem \ref{meinmaintheorem} with  
this tree $T$ we again see that Corollary \ref{meinmaintheorem+} 
is valid.   
\end{proof} 
 
\section{Conclusion} 
 
When applying ultrafilters to Ramsey theory one typically establishes  
that a set is non empty by showing that it is actually large,  
i.e. contained in a certain ultrafilter $e$.   
The Milliken-Taylor Theorem mentioned   
in the introduction states that for any  
partition $A_1,A_2,\ldots, A_r$ of $\N$ 
there exist $i$ and a sequence $\seq x$ such that 
$[\FS (\seq x)]_<^k\subset A_i.$  
In the spirit of the principle stated above, one could  
expect that after constructing the 
 first $n$ elements  $x_0,x_1, \ldots, x_{n-1}$   
of the sequence, the set of 
 possible choices of the element $x_n$ is contained 
in an ultrafilter $e$. Lemma \ref{baumramsey} gives this 
 idea a rigorous meaning. 
The combinatorial gain is that the sequence $\seq x$ can 
be forced to satisfy additional properties:   
In  
our generalizations \ref{meinmt} of  
the Milliken-Taylor Theorem the sequence may be chosen from 
a predefined IP-set in a quite general semigroup.   
In appropriate commutative semigroups  the variety of  
possible sequences is large enough to achieve the multidimensional  
extension \ref{meinmaintheorem} of the \CST. 

\section*{Acknowledgements}
The  author thanks  Vitaly Bergelson, Neil Hindman 
and Dona Strauss  for useful hints and remarks. He also thanks
the referee whose valueable
 suggestions led to a significantly more concise and elegant presentation.

\end{document}